\input amstex
\input xy
\xyoption{all}
\input epsf
\documentstyle{amsppt}
\document
\magnification=1200
\NoBlackBoxes
\nologo
\hoffset=1.5cm
\voffset=1truein

\def\Z{\bold{Z}}

\def\C{\bold{C}}

\pageheight{16cm}


\bigskip

 \centerline{\bf LOCAL ZETA FACTORS}
 
 \medskip
 
 \centerline{\bf AND GEOMETRIES UNDER $\roman{Spec}\,\bold{Z}$}
 
  \bigskip

\centerline{\bf Yuri I. Manin}

\medskip

\centerline{\it Max--Planck--Institut f\"ur Mathematik, Bonn, Germany}
  
\bigskip

{\it ABSTRACT.} The first part of this note shows that the odd  period polynomial
of  each Hecke cusp eigenform for full modular group  produces via Rodriguez--Villegas
transform ([Ro--V]) a polynomial satisfying the functional equation of zeta type
and having nontrivial zeros only in the middle line of  its critical strip. 
The second part discusses Chebyshev lambda--structure of the polynomial ring
as Borger's descent data to $\bold{F}_1$ and suggests its role in
possible relation of $\Gamma_{\bold{R}}$--factor to ``real geometry over $\bold{F}_1$ (cf. also [CoCons2]).

\bigskip

\centerline{\bf Introduction}

\medskip

In his influential seminar talk [Se], Jean--Pierre Serre stated precise conjectures about the structure
of local factors of zeta functions of algebraic varieties over arithmetic rings.
In particular, he defined the local factors at complex, resp. real, archimedean completions
of the base as multiplicative combinations of gamma functions involving Hodge numbers.
(Of course, local factors at finite primes since Weil and Grothendieck were treated in terms
of Galois representations on cohomology as characteristic polynomials of Frobenii.)

\medskip

In my seminar talks [Ma1] dedicated to the geometry and arithmetics
over Jacques Tits' mythical ``field with one element $\bold{F}_1$''  I suggested the existence of 
respective local zetas ``in characteristic one'' and noticed that Riemann's gamma--factor
at the infinite prime looks like such a local factor in characteristic one
of infinite--dimensional projective space $\bold{P}^{\infty}_{\bold{F}_1}$ 
appropriately regularized.

\smallskip

More precisely, in [Ma1] I defined
the zeta function of $\bold{P}_{F_1}^k$ as 
$$
(2\pi )^{-(k+1)}s(s-1)\dots (s-k). \eqno(0.1)
$$

\smallskip
On the other hand,  Deninger ([De]) represented  the basic  $\Gamma$--factor
at (complex) arithmetical infinity as the infinite determinant of complex Frobenius map and a regularized product
$$
\Gamma_{\bold{C}}(s)^{-1}:= \frac{(2\pi)^s}{\Gamma (s)} =\prod_{n\ge 0} \frac{s+n}{2\pi}.
\eqno(0.2)
$$
Comparing (0.1) to (0.2), I suggested that this gamma--factor, with changed sign of $s$,  might be imagined as the zeta--function of the 
infinite dimensional projective space over $\bold{F}_1$.   I did not discuss
the problem of a similar interpretation of the real gamma--factor.

\smallskip

After 1992, there was a growing body of definitions and studies of $\bold{F}_1$--geometries,
cf. surveys and comprehensive bibliography in [Lo2], [Ma2]. In particular,
Ch.~Soul\'e in [So] put on a firm ground my heuristics about local zeta factors over 
$\bold{F}_1$. In particular, natural factors of zetas of  $\bold{F}_1$--schemes 
turned out to be polynomials in $s$, satisfying a functional equation 
expressing their symmetry wrt a map $s\mapsto c-s$. In the main text,
I will use for such polynomials a generic name ``zeta polynomials'',
complementing their description by the requirement that nontrivial zeros
must lie on the vertical line at the middle of critical strip, cf. Theorem 1.3 below.

\smallskip
For other insights about $\bold{F}_1$, see [KaS], [CoCons1] and the description of  A.~Smirnov's work in [LeBr],
and about Deninger's program see [CoCons2].

\smallskip

However, the bridges between characteristics zero and one, and in particular
the  $\bold{P}^{\infty}_{\bold{F}_1}$--heuristics about (0.2)  still remain to a considerable degree
elusive.

\smallskip

In this short note, I contribute additional strokes to this mystery. 

\smallskip
In Section 1, I show
that each cusp form $f$ for $PGL(2,\bold{Z})$ which is eigenform for all Hecke operators,
besides the usual $p$--factors of its Mellin transform, produces
one more polynomial that looks like ``local zeta factor in characteristic one''.
This polynomial is obtained from the odd period polynomial of $f$ in the
same formal  way as the Hilbert polynomial
of a graded algebra is produced from its Poincar\'e series,
see [Ro--V]. Formulas (1.7) and (1.8) below suggest that this formalism
can be considered as a discrete version of the Mellin transform as well.

\smallskip

For analogies with zetas and geometric interpretations of the latter, cf. also [Go].

\smallskip

In Section 2,  I suggest how an expected gamma--bridge between characteristics zero and one could take into account
the fact that in Serre's picture gamma--factors corresponding to real and complex
infinite arithmetic primes are different. To this end, I appeal to  J.~Borger's identification of
lambda--structures on schemes with descent data to $\bold{F}_1$ ([Bo], [LeBr]),
and to the idea that Habiro rings are lifts to $\roman{Spec}\,\bold{Z}$  of ``rings of analytic functions'' in characteristic one suggested in [Ma2].
Then it turns out that two different lambda--structures on the polynomial ring, the
toric one and the Chebyshev one, faithfully reflect the difference between complex and real analytic geometry
in characteristic one.
\smallskip

Notice that lambda--structures naturally appear in several contexts, related to zetas:
see e.~g. [CoCons2], [Na] and [Ra]. It would be interesting to include Borger's philosophy
into these contexts as well.

\bigskip

\centerline{\bf 1.~ Zeta polynomials from cusp forms}

\medskip
{\bf 1.1. Period polynomials and period functions.}   Here we are considering modular forms 
with respect to  $PSL(2, \Z )$, $k$ is a positive even weight; $w:=k-2$; $S_k$ denotes the space of cusp forms, 
$M_k$ is the space all modular forms of weight $k$.

\smallskip

{\it Period polynomials} for cusp forms are defined by:
$$
r_f(z):= \int_0^{i\infty} f(\tau )(\tau - z)^{k-2}d\tau, \quad 
r_f^{\pm}(z):=\frac{r_f(z)\pm r_f(-z)}{2}.
$$
The following more general formula is valid also for Eisenstein series: if $f(z)=\sum_{n=0}^{\infty}a_ne^{2\pi ins}\in M_k$,
define its {\it Eichler integral} by
$$
\Cal{E}_f(z):=  \int_z^{i\infty}( f(\tau )-a_0) (\tau - z)^{k-2}d\tau 
= -\frac{(k-2)!}{(2\pi i)^{k-1}}\sum_{n=1}^{\infty} \frac{a_n}{n^{k-1}} e^{2\pi nz}
$$
and then define {\it period function} by
$$
r_f(z):= \Cal{E}_f(z) - z^{k-2} \Cal{E}_f(-1/z)),  \quad  r_f^{\pm}(z):=\frac{r_f(z)\pm r_f(-z)}{2}.
$$
If $f$ is not cusp form, then $r_f(z)\in z^{-1}\C [z].$

\medskip

{\bf 1.2. Spaces of period  functions.} If $g\in PSL(2,\Z )$, $g(z)=\dfrac{az+b}{cz+d}$,
the right action $|_w$ of $g$ on the space $V_w$ of polynomials $r$  of degree $\le w$ is defined by
$$
(r|_wg)(z):= (cz+d)^w r(g(z)).
$$
Let $S(z)=-1/z$, $U(z)=1-1/z$ and
$$
Y_w:= \{\,r\in V_w \, |\quad r|_w(1+S)=r |_w(1+U+U^2)=0\,\} .
\eqno(1.1)
$$
For $f\in S_k$, we have $r_f(z)\in Y_w$.
Let $Y_w^{\pm}$ mean the respective subspaces of even/odd polynomials.
\smallskip

It is well known (Eichler--Shimura) that the map $r^-:\,f\mapsto r_f^-(z)$ defines an isomorphism
$S_k\to Y_w^-$, whereas $r^+$ defines an embedding of codimension one
$S_k\to Y_w^+$.

\smallskip

Recently it was proved ([ConFaIm]) that if  $f\in S_k$ is a Hecke eigenform, then
$$
U_f(z):=\frac{r^-_f(z)}{z(z^2-4)(z^2-1/4)(z^2-1)^2}
\eqno(1.2)
$$
is a polynomial without real zeros  whose complex zeros all lie on the unit circle.
Clearly, its degree is  $e:= w-10$.

\medskip

{\bf 1.3. Theorem.}  {\it Fix an integer $d>e=w-10$ and put
$$
P_f(z):=\frac{U_f(z)}{(1-z)^d.}
\eqno(1.3)
$$
There exists a polynomial $H_f(x)\in \bold{C}[x]$ of degree $d-1$ such that
$$
P_f(z) =\sum_{n=0}^{\infty} H_f(n)z^n
$$
for $|z|< 1.$ This polynomial satisfies the functional equation
$$
H_f(x)=(-1)^{d-1}H(-d+e-x)
\eqno(1.4)
$$
and it vanishes at $x=-1,\dots , -d+e+1$. All its remaining zeros
lie on the vertical line $Re\,x=-(d-e-1)/2.$ }

\smallskip

{\bf Proof.} This is a direct application of the Proposition in sec.~3 of [Ro--V] 
(due in more general form to Popoviciu), and of its Corollary. One condition for
applicability of this Proposition is ensured by the theorem about zeros of  (1.2) from [CoFaIm].
We have only to check the functional equation (9) from this Proposition, i.~e. the identity
$$
P_f(1/z) =(-1)^dz^{d-e}P_f(z).
\eqno(1.5)
$$
Rewriting (1.5) as
$$
\frac{U_f(1/z)}{(1-1/z)^d}= (-1)^dz^{d-e}\frac{U_f(z)}{(1-z)^d}
$$
one sees that it is equivalent to
$$
U_f(1/z)=z^{-e}U_f(z)
$$
that is, in view of (1.2),
$$
\frac{r^-_f(1/z)}{z^{-1}(z^{-2}-4)(z^{-2}-1/4)(z^{-2}-1)^2} = z^{10-w} \frac{r^-_f(z)}{z(z^2-4)(z^2-1/4)(z^2-1)^2}.
\eqno(1.6)
$$
Now, from $r\,|_w(1+S)=0$ it follows that $r^-_f(1/z)=-r^-_f(-1/z)=z^{-w}r^-_f(z)$.
Inserting this into (1.6), we finally get (1.5).
\medskip
{\bf 1.4. Remarks.} a). In [ER], it was proved that all zeros of the full period polynomial of a Hecke cusp form
lie on the unit circle. Similarly, all zeros of $zr_f(z)$ for Eisenstein Hecke series lie on the unit circle.
\smallskip
However, I was unable   to fit these cases into the framework of the Rodriguez--Villegas construction,
because the analog of the functional equation (1.5) seemingly fails for the complete period polynomial.

\smallskip

b). I use the generic catchword "zeta polynomials" for polynomials of one variable satisfying
a version of functional equation such as (1.4) and ``Riemann conjecture''.  In [Ro--V], it was in particular
proved that Hilbert polynomials of certain graded rings are zeta polynomials. Golyshev ([Go])
considered rings of homogeneous functions on Fano and Calabi--Yau varieties with respect to
anticanonical or related projective embeddings and found interesting geometric
correlates of these results.

\smallskip

Moreover, comparing the formula
$$
H_f(n)=\frac{1}{2\pi i}\int_{\gamma} P_f(z)z^{-(n+1)}dz
\eqno(1.7)
$$
(where $\gamma$ is a small contour around zero) with the Mellin transform 
$$
Z_f(s)= \frac{(2\pi)^s}{\Gamma (s)} \int_0^{i\infty}f(z)\left( \frac{z}{i}\right)^{s-1} d\left(\frac{z}{i}\right)
\eqno(1.8)
$$
one sees a considerable formal analogy: morally, $H_f$ is ``discrete Mellin transform of $P_f$''.
\smallskip
In particular, the argument $n $ of $H_f$ corresponds to the classical $-s$: this is consistent  with observations
in [Ro--V] and [Go].

\smallskip

However,  finding an appropriate geometric living space for
zeta polynomials $H_f$ associated with Hecke cusp forms seemingly requires
more general realm of 
 ``geometries under $\roman{Spec}\,\bold{Z}$''.
The problem is that in most versions of  $\bold{F}_1$--geometries those zeta
polynomials that appear as zeta functions of motives over $\bold{F}_{1^n}$
have only integer zeros: cf. e.~g.~[Lo1]. To the contrary, our $H_f$ seem to come from some
non--Tate motives and geometric objects lying below $\roman{Spec}\,\bold{Z}$ but not
over  $\bold{F}_1$. I expect that they arise from the levels below $\roman{Spec}\,\bold{Z}$ 
to which such moduli stacks as
$\overline{M}_{1,n}$ can be descended. 

\smallskip

c) Notice finally that period polynomials appear also in the studies of Galois action
on the Grothendieck--Teichmueller groupoid: see [Sch], [Hai] and [Po] and references therein about their role in
Hodge realisations. One can guess that a special role of period polynomials of Hecke eigenforms 
will become clearer in the light of \'etale setting.

\medskip

\bigskip
\centerline{\bf 2. Habiro Lambda--Rings}

\medskip

{\bf 2.1. Habiro rings.}  The Habiro ring $\Cal{H}$ of one variable over $\bold{Z}$
is defined as the projective limit of quotient rings $\bold{Z} [q]/ (f(q))$ where $f(q)$ runs over the multiplicative set
of monic polynomials whose all roots are roots of unity. This ring was introduced and studied
in [Hab], and in [Ma2] it was suggested to consider it as ``{\it the ring of analytic functions
on $\bold{G}_m$ lifted from $\bold{F}_1.$}'' In fact,  $\bold{Z}[q]$ is naturally embedded into the Habiro completion $\Cal{H}$,
and $q$ becomes invertible there, so that $\Cal{H}$ can be also defined as a completion
of $\bold{Z}[q,q^{-1}]$.
One can extend this definition to the case of several invertible
variables that is, functions on tori.

\medskip

{\bf 2.2. Lambda--rings.} J.~Borger developed in [Bo] the idea to interpret  Grothendieck lambda--structures on schemes
as general descent data to $\bold{F}_1$.  It is  therefore natural to expect that the Habiro ring admits a natural lambda--structure.\smallskip
Here we will be concerned only with commutative rings $A$
flat over $\bold{Z}$ in which case a lambda--structure can be considered simply
as a system of commuting lifts of Frobenii: ring homomorphisms $\psi^p:\,A\to A$
for each prime $p$ such that $\psi^p(x) \equiv x^p\, \roman{mod}\, pA$ for all $x\in A$ and $\psi^{p_1}\psi^{p_2}=\psi^{p_2}\psi^{p_1}$. In particular, we can define $\psi^k:\,A\to A$ for all positive integers $k$
by multiplicativity.
\smallskip

The most natural lambda--structure on  $\bold{Z} [q]$ and $\bold{Z}[q,q^{-1}]$  is determined by $\psi^k(q)=q^k,$ and since it is compatible
with the projective limit over the system of cyclotomic polynomials in $q$,
it is inherited by the Habiro ring. We will call this structure {\it toric one}.

\smallskip

However, the polynomial ring $\bold{Z} [r]$ admits one more lambda--structure discovered by Clauwens ([Cl]).
In this structure, 
$$
\psi^k (r ):=T_k ( r )
$$
where $T_k$ is the $k$--th Chebyshev polynomial. Our next result describes a subring $\Cal{H}_0\subset \Cal{H}$
which can be endowed with Chebyshev lambda--structure.

\medskip

{\bf 2.3. Proposition.} {\it (i) Consider in the Habiro ring $\Cal{H}$ the subring $\Cal{H}_0$
defined as the completion of the polynomial subring $\bold{Z}[r]$, where
$$
 \quad r:=1+q+\sum_{n=1}^{\infty} q^n\cdot (1-q)\dots (1-q^n).
 \eqno(2.1)
$$
This subring is invariant with respect to the standard lambda--structure $\psi^k$, which induces on
this subring, in terms of the coordinate $r$, the Chebyshev lambda--structure.
\smallskip

 (ii) $\Cal{H}_0$ is strictly smaller than $\Cal{H}$.}

\medskip

{\bf Proof.} (i) In $\Cal{H}$, we have the convergent expression for $q^{-1}$ (see [Hab], Prop. 7.1):
$$
q^{-1}=1+\sum_{n=1}^{\infty} q^n\cdot (1-q)\dots (1-q^n)
$$
Hence $r=q+q^{-1}$. Moreover, using one of the definitions of Chebyshev polynomials, we see that
$$
\psi^k ( r )=q^k+q^{-k} = T_k(q+q^{-1}) =  T_k ( r).
$$
\smallskip

(ii) In order to see that  $\Cal{H}_0$ is strictly smaller than $\Cal{H}$, we can use the following result due
to Habiro. Any element of $\Cal{H}$ determines a function on the set of roots of unity $\mu_{\infty}$
with values in $\bold{Z}[\mu_{\infty}]$, and the resulting map
$$
\Cal{H}\to \roman{Map} (\mu_{\infty}, \bold{Z}[\mu_{\infty}])
$$
is an embedding (see [Hab]). The element $q$ corresponds to the tautological map
$\mu_{\infty}\to \bold{Z}[\mu_{\infty}]$.

\smallskip

Then  all elements of $\bold{Z}[r]$  become functions invariant under the involution
$\zeta\to\zeta^{-1}$ of $\mu_{\infty}$ and their values  are as well invariant.
This property holds after the completion. Hence $q\notin \Cal{H}_0$.

\smallskip

Notice that for each complex embedding of $\mu_{\infty}$ and any $\eta \in  \mu_{\infty}$,
$\eta+\eta^{-1}$ is real. This is why we referred to ``real analytic geometry over $\bold{F}_1$.''

\bigskip

\centerline{\bf References}
\medskip

[Bo] J.~Borger. {\it $\Lambda$--rings and the field with one element.} arXiv:0906.3146

\smallskip

[Cl] F. J. B. J. Clauwens. {\it Commuting polynomials and  $\lambda$--ring structures on
$\bold{Z}[x]$.} J.~Pure Appl.~ Algebra, 95(3), 1994,  261--269.

\smallskip

[CoCons1] A.~Connes, C.~Consani. {\it  Schemes over $\bold{F}_1$ and zeta functions.} Compos. Math. 146, no. 6, 2010, 1383--1415.

\smallskip

[CoCons2] A.~Connes, C.~Consani. {\it Cyclic homology, Serre's local factors and the $\lambda$--operations.}
arXiv:1211.4239

\smallskip

[ConFaIm] J.~B.~Conrey, D.~W.~Farmer, \"O.~Imamoglu. {\it The nontrivial zeros of period
polynomials lie on the unit circle.} Int.~Math.~Res.~Not.~, no. 20, 2013, 4758--4771.
arXiv:1201.2322
\smallskip

[De] C.~Deninger. {\it On the $\Gamma$--factors attached to motives.} Inv.~Math. 104, 1991, 245--261.

\smallskip
[ER] A.~El--Guindy, W.~Raji. {\it Unimodularity of zeros of period polynomials
of Hecke eigenforms.} Bull.~Lond.~Math.~Soc. 46, 2014, 528--536.
\smallskip
[Go] V.~Golyshev. {\it The canonical strip, I.} arXiv:0903.2076

\smallskip

[Hab] K.~Habiro. {\it Cyclotomic completions of polynomial rings.} Publ.~RIMS,
Kyoto Univ., 40, 2004, 1127--1146.

\smallskip

[Hai] R.~Hain. {\it The Hodge--de Rham theory of modular groups.} arXiv:1403.6443

\smallskip

[KaS]. M.~Kapranov, A.~Smirnov. {\it Cohomology determinants and reciprocity laws:
number field case.}  Unpublished manuscript, 1996, 15 pp.
\smallskip

[LeBr] L.~Le Bruyn.  {\it Absolute geometry and Habiro topology.} arXiv:1304.6532

\smallskip
[Lo1] O.~Lorscheid. {\it Functional equations for zeta functions of $\bold{F}_1$--schemes.}
C.~R.~Ac~Sci.~Paris, 348 (21--22), 2010, 1143--1146. arXiv 1010.1754
\smallskip

[Lo2] O.~Lorscheid. {\it A blueprinted view on $\bold{F}_1$--geometry.}
arXiv:1301.0083

\smallskip

[Ma1] Yu.~Manin. {\it  Lectures on zeta functions and motives (according to Deninger and Kurokawa).}
 In: Columbia University Number Theory Seminar (1992),  Ast\'erisque 228,
1995, 121--164.

\smallskip
[Ma2] Yu.~Manin. {\it  Cyclotomy and analytic geometry over $F_1$.} In: Quanta of Maths. Conference in honour of Alain Connes.
Clay Math. Proceedings, vol. 11, 2010, 385--408.
 Preprint math.AG/0809.2716. 28 pp.
 
 \smallskip
 
 [Na] N.~Naumann. {\it Algebraic independence in the Grothendieck ring of varieties.}
 Trans. AMS, 359(4): 1652--1683 (electronic), 2007.

\smallskip

[Po] A.~Pollack. {\it Relations between derivations arising from modular forms.}
http://dukespace.lib.duke.edu/dspace/handle/10161/1281
\smallskip

[Ra] N.~Ramachandran. {\it Zeta functions, Grothendieck groups, and the Witt ring.}
arXiv:1407.1813

\smallskip

[Ro--V] F.~Rodriguez--Villegas. {\it On the zeros of certain polynomials.}
Proc.~AMS, Vol.~130, No 8, 2002, 2251--2254.

\smallskip

[Se] J.-P.~Serre. {\it Facteurs locaux des fonctions z\^eta des vari\'et\'es alg\'ebriques
(d\'efinitions et conjectures).} S\'eminaire Delange--Pisot--Poitou 11, 1969--1970, exp. 19, 1--15.

\smallskip

[Sch] L.~Schneps. {\it On the Poisson bracket on the free Lie algebra in two generators.}
J.~Lie Theory 16, no. 1, 2006, 19--37. 

\smallskip
[So] Ch.~Soul\'e. {\it Les vari\'et\'es sur le corps \`a un \'el\'ement.} Mosc. Math. J., 4(1), 2004, 217--244.


\enddocument